\documentclass[a4paper,11pt]{article}
\usepackage{times}
\usepackage{graphicx}
\usepackage{amsfonts}
\usepackage{amsmath, amsthm, amssymb}
\usepackage{dsfont}
\usepackage{ucs}
\usepackage[utf8x]{inputenc}
\usepackage{color}
\usepackage{geometry}
\usepackage{hyperref}
\usepackage{ifthen}
\title{On an explicit representation of the solution of linear stochastic partial differential equations with delays}
\author{Mathieu Galtier\footnote{corresponding author: mathieu.galtier@inria.fr} \and Jonathan Touboul}
\date{}

\def \N {\mathbb{N}}

\def \F {\mathcal{F}}

\def \U {\mathcal{U}}
\def \E {\mathcal{E}}
\def \V {\mathcal{V}}
\def \T {\mathcal{T}}
\def \Lop {\mathcal{L}}
\def \Co {\mathcal{C}}

\def \X {\mathcal{X}}
\def \B {\mathcal{B}}
\def \D {\mathcal{D}}
\def \R {\mathbb{R}}
\def \C {\mathcal{C}}
\def \eqdef {\stackrel{def}{=}}

\newtheorem{thm}{Theorem}[section]
\newtheorem{prop}[thm]{Property}

\begin{document}
\maketitle
%\vspace{-1cm}
\begin{center}
 NeuroMathComp team, INRIA / ENS, 23 avenue d'Italie, 75013 Paris, France
\end{center}
\begin{center}
September 2011
\end{center}

\begin{abstract}
Based on the analysis of a certain class of linear operators on a Banach space, we provide a closed form expression for the solutions of certain linear partial differential equations with non-autonomous input, time delays and stochastic terms, which takes the form of an infinite series expansion.

\medskip

\end{abstract}

\noindent{\bf Introduction}\\
Linear differential systems are ubiquitous in pure and applied mathematics, either as models, approximations, but also because the stability of solutions of nonlinear differential systems reduces to the study of linear systems. Such systems might include stochastic terms (see \cite{mao:08}), temporal delays (see \cite{hale-lunel:93}), and also encompass the case of partial differential equations. Apart from the simplest linear finite-dimensional differential equations, finding closed forms expressions for the solutions of general linear differential systems is very complex. In this paper, based on the treatment of evolution equations as algebraic equations in a suitable Banach space, we propose a closed-form expression for the solution of linear, non-autonomous, stochastic, time-delayed partial differential systems. Application of this framework to several classical examples such as the delayed Ornstein-Uhlenbeck process or the stochastic heat equation are developed in sections~\ref{sec:examples} and \ref{sec:heat equation}. This expression is especially useful to understand the dynamics of weakly connected linear learning neural networks, problem which motivated the development of this more general framework (sections \ref{sec: motivation} and \ref{sec:WCN}). 

%\vspace{-0.6cm}
\section{{Motivation: neural networks with learning}}
\label{sec: motivation}
%\vspace{-0.5cm}
{
Consider a neural network made of $n\in \N^*$ neurons described by their membrane potential at time $t$ written $V(t) \in \R^n$. The connectivity matrix $W\in \R^{n \times n}$ contains the weights of the connections between each pair of neurons. Learning in biological networks corresponds to the slow evolution of the connectivity. Neural network with Hebbian learning can be modeled by the very simple coupled model:
%\vspace{-0.2cm}
\begin{equation}
\left\{
\begin{array}{ll}
 \epsilon \dot{V} & = -l V + W \cdot V + I(t)\\
 \dot{W} &= - \kappa W + V \otimes V 
\end{array}
\right.
\label{eq: linear learning neural network}
\end{equation}}

%\vspace{-0.3cm}
\noindent {where the positive parameters $\epsilon$ governs the ratio between the potential and the connectivity time scales, $l$ the leak conductance and $\kappa$ the learning ratio, and $I(t) \in \R^n$ is a $\tau$-periodic external input. The symbol $\otimes$ denotes the tensor product: $\{V\otimes V\}_{ij} = V_i V_j$. A companion paper (in preparation) proves the well-posedness, and analyzes the dynamics and the equilibrium connectivity of such systems in the limit $\epsilon \to 0$, i.e. when learning occurs on a much slower scale than the activity. In that asymptotic regime, the connectivity is well approximated by the equation $\dot{\bar{W}} = - \kappa \bar W + \frac{1}{\tau}\int_0^\tau V_{\bar W}(t) \otimes V_{\bar W}(t) dt$ where $V_{\bar W}$ is the periodic solution of $\dot V  = -l V + \bar W \cdot V + I(t)$ where $\bar W$ is assumed constant. Seeing $V_{\bar W} \in C^{1}([0,\tau[,\R^n)$ as a ``semi-continuous'' matrix of size $n \times [0,\tau[$, such that ${V_{\bar W}}_{it} = ({V_{\bar W}}(t))_i$, leads to the formulation
%\vspace{-0.2cm}
$$
 \dot{\bar W} = - \kappa \bar W + \frac{1}{\tau} V_{\bar W}\cdot V_{\bar W}'
$$}
%\vspace{-0.4cm}

{The characterization of the solutions of these equations is therefore important to understand what is learned by such networks. The following results will provide a way to compute the equilibria as shown in section~\ref{sec:WCN}. }
% Thus having a matrix-like representation of $V_{\bar W}$ would make it possible to solve the previous equation. In the following we derive such a representation which takes the form of a weakly-connected expansion: we are able to express $V_{\bar W}$ as an infinite converging sum whose terms are vanishing even more when the connectivity is small (which is a biologically plausible assumption). Indeed, theorem \ref{thm: weakly expansion} makes it possible to compute explicitely the equilibrium connectivity of system \eqref{eq: linear learning neural network}.}

%\vspace{-0.6cm}
\section{Framework and General Result}
\label{sec:theoretical}
%\vspace{-0.5cm}
The framework we develop here is based on extending notions of matrix calculus to infinite dimensional spaces. The linearity of the equation motivates to extend some finite-dimensional linear algebra and matrix concepts to infinite-dimensional spaces.

We consider in the manuscript linear equations in a Banach space $\C$ of real functions of time $t$ and a variable $x\in E$, called \emph{space} variable, where $E$ can either be a finite set $\{1,\ldots,N\}$ (in which case $\C$ is equivalent to the space of $\R^N$-valued functions), countable or continuous, typically $\R$, in which case $\C$ is a space of two-variables functions. The particular problem under consideration governs the choice of the space $\C$, in particular including regularity or integrability properties (typically $\C$ is a $L^p$ or a Sobolev space). Similarly to a matrix notation, we denote the value of $X\in \C$ at $(x,t) \in E \times \R$ by $X_{xt}$. 

Let $\E$ denote the space of bounded linear operators on $\C$. We are interested in solving equations of type $\Lop X = B$ where $\Lop \in \E$ (this operator may involve differentials in time and/or space) and $B\in \C$. We will restrict the study to a class of operators of a particular form we now detail. To this end, we introduce two kinds of linear operators on $\C$: the \emph{space} operators $L$ acting on the first (space) variable, i.e. linear operators on $\R^E$. If $E$ is finite, this set is reduced to the matrices. If $E$ is equal to $\R^d$, it contains all the linear operators acting on functions of the space variable, in particular, under suitable regularity conditions, integral or differential operators.  The action of the space operators $L$ on a function $X\in \C$ is denoted $L\cdot X$ (acting on the left). The \emph{time} operators essentially act on the second (time) variable, and the transform might depend on the space variable $x$. In other words, these transforms $\mathcal{R}$ can be represented by a family of operators $(\mathcal{R}_x, x\in E)$ such that for any $x$, $\mathcal{R}_x$ is a linear operator on $L^{2}(\R)$. The action of a time operator $\mathcal{R}$ on $X \in \C$ is written $X \cdot \mathcal{R}$ (acting on the right). In the paper, we will mainly be interested in diagonalizable time operators. Diagonal operators in the time domain are operators $\mathcal{R}$ whose action can be written in the form $(X\cdot \mathcal{R})_{xt}=r(x,t) X_{x,t}$. This class includes for instance all linear differential time operators, which are diagonalizable in the Fourier basis. Another class of time operators we will be considering is the class $\Co$ of convolution operators with respect to time. Given a finite measure $g$ of $\R$, the convolution operator $\T_g \in \Co$ associated with $g$ is defined as $\big(X \cdot \T_g \big)_{xt} = \int_{-\infty}^{\infty} X_{x (t-s)} dg(s)$. Such operators are generalizations of Toeplitz matrices generated by $g$, with, loosely speaking, infinitely many rows and columns. An important property of the convolution operators is that they are diagonal in the Fourier basis.

For $L$ a space operator and $\mathcal{R}$ a time operator, we define the Kronecker product $L \otimes \mathcal{R}$ as the mixed operator of $\E$ such that $\big( L \otimes \mathcal{R} \big)(X) = L\cdot (X \cdot \mathcal{R})$. Note that the product becomes associative when $\mathcal{R}$ is a convolution operator which will be the case in section \ref{sec:applications}. This definition extends the property of vectorization of the Kronecker product of matrices in linear algebra (see e.g.~\cite{brewer:78}).

\medskip

The main technical result of the paper is given in the following:

\begin{prop}\label{thm: Kronecker inversion}
  Let $\Lop=A \otimes \B + Id_{\C} \otimes \D$ be a linear operator, for $A$ a space {operator} and $\B,\,\D$ co-diagonalizable time operators, with $\B$ invertible. For the sake of simplicity, we assume that they are diagonal in the natural time basis, and denote for $x\in E$, $\B_x = diag_{t\in \R}\Big(b(x,t)\Big)$ and $\D_x= diag_{t\in \R}\Big(d(x,t)\Big)$. We assume that $\inf_{x,t}\vert b(x,t) \vert >0$ and the spectral condition:
\begin{equation}
\label{eq: spectral condition} 
	\exists\, l \in \R^* \text{ such that } \; \lambda \eqdef \frac{\|W\|}{\inf_{x,t}\Big\{\big|l - \frac{d(x,t)}{b(x,t)} \big| \Big\}} < 1 ,
\end{equation}
where $W \eqdef l \,\textrm{Id}_{\R^E} + A$ and $\|W\| = \sup_{X\neq 0} \frac{\|W \cdot X\|}{\|X\|}$ is the operator norm.
Then $A \otimes \B + Id_{\C} \otimes \D$ is invertible and its inverse reads:
\begin{equation}\label{eq:Inverse}
 \Big(A \otimes \B + Id_{\C} \otimes \D\Big)^{-1} = -\sum_{k = 0}^{+\infty} W^k \otimes diag_{t\in \R}\Big(\frac{1}{b(x,t)(l-\frac{d(x,t)}{b(x,t)})^{k+1}}\Big).
\end{equation}
% where $\mathcal{U}=\F \cdot diag_{x\in E} \Big(\frac{l}{l b(x) - d(x)}\Big) \cdot \F^{-1}$ and $\mathcal{V}=\F \cdot diag_{x\in E}\Big(\frac{l b(x)}{l b(x) - d(x)}\Big)\cdot \F^{-1}$
\end{prop}
	{\noindent {\bf Remark} {\it The spectral condition is merely a technical sufficient condition for the convergence of the series. The relatively formal setting and assumptions will become clearer in the applications, section~\ref{sec:applications}.}}

\begin{proof}
 The direct introduction of the inverse can appear artificial at first sight. However, this formula is a natural extension of the discrete-time case where direct linear algebra and Kronecker products calculations quite simply provide a closely related expression the interested reader can readily derive. 

In order to prove the proposition, we first need to prove that the operator indeed exists, and that it constitutes the inverse of $\Lop$. {It is easy to show that under the assumption of the proposition that the sequence of operators in $\E$ defined by: $M_N \eqdef -\sum_{k = 0}^{N} W^k \otimes diag_{t\in \R}\Big(\frac{1}{b(x,t)(l-\frac{d(x,t)}{b(x,t)})^{k+1}}\Big)$ constitutes a Cauchy sequence in $\E$. Since $\C$ is a Banach space, so is $\E$, and hence the sequence $(M_n)_n$ converges. The limit of this sequence is our inverse candidate, and is denoted as the infinite series \eqref{eq:Inverse}.}

In order to prove that this limit is indeed the inverse of $\Lop$, we compute the limit of $(M_N \circ \Lop) X$ (or similarly $(\Lop \circ M_N) X$) for a given $X \in \C$. It is easy to show, developing the series, that we have:
\[((M_N \circ \Lop) X)_{xt} = - \sum_{k=0}^N -W^k \cdot \frac{X_{.,t}}{\Big(l-\frac{d(.,t)}{b(.,t)}\Big)^k} + W^{k+1} \cdot \frac{X_{.,t}}{\Big(l-\frac{d(.,t)}{b(.,t)}\Big)^{k+1}} = X_{xt} - W^{N+1}\cdot \frac{X_{.,t}}{\Big(l-\frac{d(.,t)}{b(.,t)}\Big)^{N+1}}\]
{where $Y_{.,t}$ for $Y \in \C$ denotes the application $E\mapsto \R$ such that $Y_{.,t}(x)=Y_{xt}$.
Here again, the assumptions of the proposition ensures that the second term vanishes as $N$ goes to infinity. }
\end{proof}

%\vspace{-0.6cm}
\section{Application to solving linear time-delayed Stochastic Partial Differential Equations}\label{sec:applications}
%\vspace{-0.5cm}
In this section we make explicit the use of the inversion formula~\eqref{eq:Inverse} in the case of linear delayed, stochastic, partial differential equations. Several examples with different convolution operators will illustrate the main result of the section stated in theorem~\ref{thm: weakly expansion}. 
%\vspace{-0.4cm}
\subsection{General Result}
%\vspace{-0.4cm}
Let $\X$ be a Hilbert space, typically $\R^n$ for $n\in \N$, $L^2(\R^n)$ or a Sobolev space of applications on $\R^n$. We consider a probability space $(\Omega, \mathcal{F}, \mathbb{P})$ satisfying the usual conditions and $B$ a standard adapted $\X$-Brownian motion (for the existence and properties of this object in infinite-dimensional spaces, see~\cite[Chapter 4]{da-prato:92}). We aim at solving the non-autonomous time-delayed stochastic differential equation:
\begin{equation}
\left\{
\begin{array}{c}
dX = \big(A\cdot (X \ast g) + I\big)\ dt + \Sigma\cdot dB\\
X_{|_{\R_-}} = \zeta_0 \in L^2_{\X}(\R_-^*)
\end{array}
\right.
\label{eq: differential SDE}
\end{equation}
with $\Sigma : \X\mapsto \X$ linear,  $I \in C(\R_+,\X)$ an external input, $g$ a finite measure of the real line supported on $\R_+$, i.e. a causal measure, and $\ast $ denoting the convolution. Existence and uniqueness of weak solutions for such equations is ensured, see e.g.~\cite{mao:08,da-prato:92}. We consider the case where the system has a unique strong solution. In the case where $\X=\R^n$, this occurs under the assumptions of the section, see e.g.~\cite[Chapter 5]{mao:08}, and in the infinite-dimensional case, we need to assume that $B$ is a genuine Wiener process (i.e. the trace of the covariance matrix is finite, and the initial condition is in the domain of A, see~\cite{da-prato:92}). The solution of this stochastic differential equation at time $t\in \R_+$ is defined by the integral equation $X(t) = \zeta_0(0) + \int_0^t \big(A\cdot (X\ast g)(s) + I(s)\big)\ ds + \int_0^t \Sigma\cdot dB$ and $X_{|_{\R_-}} = \zeta_0$.

This problem can be set in the framework described in section~\ref{sec:theoretical} using {a transformation inspired by the classical} Fourier transform of the solution in the time domain. To perform this transformation {rigorously in our particular stochastic setting, we} stop our processes at a finite time $\tau>0$. We define $X_\tau: t \in \R \rightarrow \mathds{1}_{[0,\tau]}(t) X(t)$ the restriction of $X$ to the compact support $[0,\tau]$ and null elsewhere. Similarly, define $I_\tau = \mathds{1}_{[0,\tau]} I$ and $dB_\tau = \mathds{1}_{[0,\tau]} dB$. We have:
%\vspace{0.2cm}
\begin{thm}
\label{thm: weakly expansion}
For all $\tau \in \R_+$, choose $l \in \R^*$ and $W$ a space operator such that $W = l\,Id_{\C} + A$. If the spectral condition \eqref{eq: spectral condition} is satisfied, i.e. in the present case {$\|W\| < \inf_\xi\{\big|l + \frac{2i\pi \xi}{\hat{g}(\xi)}\big|\}$}, then the solution of equation (\ref{eq: differential SDE}) is given by
\begin{equation}
X_\tau =  \sum_{k=0}^{+\infty} W^k\cdot \Big(\zeta_0(0) \delta_0 +  \widetilde{I}_\tau + \Sigma\cdot dB_\tau\Big)\cdot \U\cdot \V^k
\label{eq: weakly connected expansion}
\end{equation}
where $\mathcal{U}=\F\cdot diag_{\xi\in \R}\Big(\frac{1}{l \hat{g}(\xi) +2i\pi \xi}\Big)\cdot \F^{-1}$, $\mathcal{V}=\F\cdot diag_{\xi\in \R}\Big(\frac{\hat{g}(\xi)}{l \hat{g}(\xi) +2i\pi \xi}\Big)\cdot \F^{-1}$ and $\widetilde{I}_\tau = I_\tau + A \cdot (\zeta_0 \ast g)$. The notation $\big(dB_\tau \cdot \U\big)_{xt}$ stands for the stochastic integral $\int_0^{min(t,\tau)} \U_{st} dB(s)$ which is square integrable on $[0,\tau[$.
\end{thm}
%\vspace{0.2cm}
% {There is no clear explanation of the meaning of the spectral condition, yet it could be understood as a requirement of non pathological intereferences between the oscillations coming from the spatial and temporal operators.}
\noindent{{{\bf Remark:}{\it The convergence of the series~\eqref{eq: weakly connected expansion} occurs as soon as the spectral condition~\eqref{eq: spectral condition} is satisfied on the subspace spanned by $W^k\cdot (\zeta_0(0) \delta_0 +  \widetilde{I}_\tau + \Sigma\cdot dB_\tau)\cdot \U\cdot \V^k$.}}}\\
\begin{proof}First, note that $A\cdot (X \ast g) = A\cdot (X_\tau \ast g) + A\cdot (\zeta_0 \ast g)$ yielding the equation on $X_\tau$: $dX_\tau = \big(A\cdot (X_\tau \ast g) + \widetilde{I}_\tau\big)\ dt + \Sigma\cdot dB_\tau$. Thus, the initial condition on $X$ acts as an external input on $X_\tau$. In the deterministic finite-dimensional case, it is well known that differential operators are diagonal in the Fourier basis. Based on this result, we introduce the Fourier transform $\F$ of equation \eqref{eq: differential SDE} for a fixed $\omega \in \Omega$. As mentioned, for almost all $\omega \in \Omega$, the processes involved are bounded, hence the function of time, on the compact interval $[0,\tau]$, is square integrable in time. Let $Z^\xi: t \in \R \rightarrow e^{-2i\pi t\xi} X_\tau(t)$ for $\xi\in \R$ the Fourier variable and $X$ is the unique solution of equation \ref{eq: differential SDE}. It\^o formula yields {for $t < \tau$}
\begin{equation}
dZ^\xi(t) = \Big(-2i\pi \xi Z^\xi(t) + A\cdot \big(Z^\xi\ast g\big)(t)+e^{-2i\pi t \xi} \widetilde{I}_\tau(t)\Big)dt + e^{-2i\pi t \xi} \Sigma\cdot dB_\tau(t)
\label{eq: after Ito}
\end{equation}
% which reads, for all $\tau \in \R_+$,
% \begin{equation}
%  Z^\xi(\tau) = X_0(0) -2i\pi \xi \int_0^\tau Z^\xi(s) ds + A.\int_0^\tau (Z^\xi \ast g)(s) ds+ \int_0^\tau e^{-2i\pi t \xi} I(s) ds +  \int_0^\tau e^{-2i\pi s \xi} dB(s)
% \end{equation}
Let us denote by $\hat{X}_\tau: \xi \in \R \to \int_0^\tau Z^\xi(s) ds$ the Fourier transform of $X_{\tau}$ {and $ \hat{I}_\tau:\xi \to \int_0^\tau e^{-2i\pi t \xi} \tilde{I}_\tau(t) dt$}. The process $\hat{B}_\tau$ is the well-defined stochastic integral $\int_0^\tau e^{-2i\pi t \xi} dB(t)$. The integral form of equation \eqref{eq: after Ito}, using the fact that the convolution is diagonal in the Fourier basis, denoting $\hat{D} = diag_{\xi\in\R}\big(-2i \pi \xi\big)$ and $\hat{\mathcal{G}} = diag_{\xi\in\R}\big(\hat{g}(\xi)\big)$, leads to the functional equation:
\begin{equation*}
Z^{\cdot}(\tau) - Z^{\cdot}(0) = A\cdot \hat{X}_\tau\cdot \hat{\mathcal{G}} + \hat{X}_\tau\cdot \hat{D} + \hat{I}_\tau + \Sigma\cdot \hat{B}_\tau
% \label{eq: sylvester sto}
\end{equation*}
Applying proposition \ref{thm: Kronecker inversion} for a fixed $\omega \in \Omega$ where $\C$ is the set of square integrable functions on $[0,\tau[$ which is a Banach space, we obtain:
\begin{equation}
\hat{X}_\tau = \sum_{k=0}^{+\infty} W^k\cdot \Big(- Z^{\cdot}(\tau) + Z^{\cdot}(0) + \hat{I}_\tau + \Sigma\cdot \hat{B}_\tau \Big)\cdot diag_{\xi\in \R}\Big(\frac{1}{\hat{g}(\xi)(l+\frac{2i\pi \xi}{\hat{g}(\xi)})^{k+1}}\Big)
\label{eq: after inversion Kronecker}
\end{equation}
We now take the inverse Fourier transform of this expression by applying the time operator $\F^{-1}$. First of all, we perform the inversion on the terms $\hat{I}_\tau = \widetilde{I}_\tau\cdot\F$. It is easy to show that $\hat{I}_\tau \cdot diag\big(\frac{1}{\hat{g}(\xi)(l+\frac{2i\pi \xi}{\hat{g}(\xi)})^{k+1}}\big)\cdot \F^{-1} = \widetilde{I}_\tau\cdot \U\cdot \V^{k}$. Similarly, the term $\Big(\hat{B}_\tau\cdot diag\big(\frac{1}{\hat{g}(\xi)(l+\frac{2i\pi \xi}{\hat{g}(\xi)})^{k+1}}\big) \cdot \F^{-1}\Big)_{.t}$ can be written $ dB_\tau \cdot \U\cdot \V^k$.

Moreover, for $x \in \{0,\tau\}$ an easy computation shows that {$\Big(Z^{\cdot}(x)\cdot diag\big(\frac{1}{\hat{g}(\xi)(l+\frac{2i\pi \xi}{\hat{g}(\xi)})^{k+1}}\big)\cdot \F^{-1}\Big)_{.t} = \big(X(x)\delta_x\big)\cdot \U\cdot \V^k$}.
% \begin{equation}
%  \Big(Z(x).diag\big(\frac{l^{k+1}}{\hat{g}(x)(l+\frac{2i\pi x}{\hat{g}(x)})^{k+1}}\big).\F^{-1}\Big)_{.t}
% = \int_{-\infty}^{+\infty} e^{2i \pi \xi (t-x)} \frac{l^{k+1}}{\hat{g}(\xi)(l+\frac{2i\pi \xi}{\hat{g}(\xi)})^{k+1}} X(x)\ d\xi \\
% = \int_{-\infty}^{+\infty} \F^{-1}\Big(\frac{l^{k+1}}{\hat{g}(\xi)(l+\frac{2i\pi \xi}{\hat{g}(\xi)})^{k+1}} \Big)(t-y)  \delta_x(y) X(y)\ dy
% = \big(X(x)\delta_x\big).\U.\V^k
% \end{equation}
Furthermore, the operators $\U$ and $\V$ are causal, i.e. if {$Y$ has a support $\subset [c, +\infty[$ then $Y\cdot \U\cdot \V^k$} also has a support $\subset [c,+\infty[$. {Indeed, $\hat{u}:\xi \mapsto \frac{1}{l\hat{g}(\xi) +2i\pi \xi}$ corresponds to the transfer function of a closed loop filter shown on the right, and hence $\U$ is clearly causal since $g$ is. }

\hspace{-0.4cm}\begin{minipage}{0.70\textwidth}
{$\V$ is also causal as the convolution of $g$ and $\U$. This implies that the contribution of $Z(\tau)$ vanishes in equation (\ref{eq: after inversion Kronecker}) since it has its support in $[\tau, \infty]$. }
\end{minipage}\quad
\begin{minipage}{0.25\textwidth}
 \includegraphics[width=.8\textwidth]{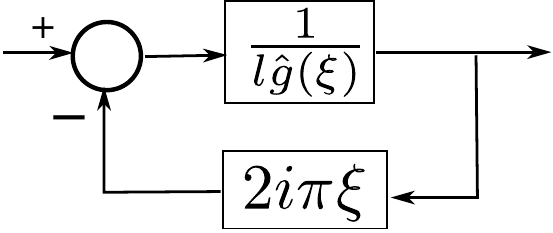}
\label{fig: closed loop}
 % closed_loop.pdf: 197x142 pixel, 72dpi, 6.95x5.01 cm, bb=0 0 197 142
\end{minipage}

%  \begin{figure}[ht]
%  \centering
%  \includegraphics[width=0.2\textwidth]{../closed_loop.pdf}
%  % closed_loop.pdf: 197x142 pixel, 72dpi, 6.95x5.01 cm, bb=0 0 197 142
%  \label{fig: closed loop}
% \end{figure}

% 
% \begin{equation}
%  \Big(\hat{dB}_\tau.diag\big(\frac{l^{k+1}}{\hat{g}(x)(l+\frac{2i\pi x}{\hat{g}(x)})^{k+1}}\big).\F^{-1}\Big)_{.t}
% = \int_{-\infty}^{+\infty} e^{2i \pi t \xi} \frac{l^{k+1}}{\hat{g}(\xi)(l+\frac{2i\pi \xi}{\hat{g}(\xi)})^{k+1}} \int_0^{\tau} e^{-2i\pi \xi s}dB(s)\ d\xi\\
% =  \int_0^{\tau} \F^{-1}\Big(\xi \mapsto \frac{l^{k+1}}{\hat{g}(\xi)(l+\frac{2i\pi \xi}{\hat{g}(\xi)})^{k+1}}\Big)(t-s)\ dB(s) =  dB_\tau.\U.\V^k
% \end{equation}
\end{proof}
\vspace{-1cm}
\subsection{Computational Remarks}
%\vspace{-0.4cm}

Truncations of the formula \eqref{eq: weakly connected expansion} provides approximations of the solution of system \eqref{eq: differential SDE}. We observe that the smaller the difference between the spatial operator $A$ and a multiple of the identity, the more accurate a truncation of this expansion. In other words, this expansion is particularly useful if $A$ is a small perturbation of the (scaled) identity (e.g. the case of weakly connected linear neural networks).

This representation allows development of new numerical schemes for the simulations of the solutions of system \eqref{eq: differential SDE}. For simplicity, consider the case where $E= \{1,\cdots,n\}$. To approximate the solution over the interval $[0,\tau]$ define a time step $\Delta t$ an number of points $T= \tau/\Delta t \in \N$ and replace $\U$ and $\V$ by the Toeplitz square matrices $\tilde{\U}$ and $\tilde{\V}$, generated by $i \in \{0,\cdots, T-1\}  \mapsto \int_{i \Delta t}^{(i+1)\Delta t} u(s) ds$, where $u$ is the function generating $\U$ (and similarly for $\V$). The number of operations needed is $\mathcal{O}((k+1)nT(n + \ln T))$ since the product with a Toeplitz matrix, as a convolution, has a cost $\mathcal{O}(T \ln T)$. This scheme has a first order accuracy, $\mathcal{O}(dt^\gamma + dx + \lambda^{k+1} )$ where $\gamma$ is equal to $1$ for deterministic equations or $\frac{1}{2}$ if stochastic. {In comparison, the Euler-Maruyama method has a complexity of $\mathcal{O}(T\big(n^2 + n\,\frac{\theta}{dt} ln(\frac{\theta}{dt})\big) )$ where $\theta$ is the support of $g$ and an accuracy of $\mathcal{O}\big(dt^\gamma + dx\big)$, comparable to the expansion method in both aspects.}

{Two interesting advantages of the expansion over Euler-like methods are that (i) it is parallelizable and (ii) it appears to be numerically very stable, i.e. large $\Delta t$ do not lead to a diverging scheme.}

%\vspace{-0.4cm}
\subsection{Examples}
%\vspace{-0.4cm}
\label{sec:examples}
Let us now treat some classical problems that are solved in the present framework.
\begin{itemize}
 \item \textbf{Ornstein-Uhlenbeck process:} The simplest example is the Ornstein-Uhlenbeck process with no delays (i.e. $g = \delta_0$). In that case, $\hat{g} = 1$, and therefore, for all $l \in \R$, $\inf_\xi\{|l +2i\pi \xi |\} = |l|$ and the expansion is valid if there exists $l\in \R^*$ such that $\|l + A\|<\vert l \vert$, i.e. for any operator $A$ whose spectrum is bounded {and entirely contained in the left or right half plane}. {For negative matrices $A$ (i.e. $l>0$)} $\T_h=\U = \V$ is a Toeplitz operator generated by the function $h:z \rightarrow e^{-lz}H(z)$ with $H$ the Heaviside function. Therefore, the solution of $\dot{X} = A\cdot X + I + \Sigma\cdot dB$ can be written as
$
X_\tau = \sum_{k=0}^{+\infty} W^k\cdot (X_0 \delta_0 + I_\tau + \Sigma\cdot dB_\tau )\cdot \T^{k+1}_h.
$
If $A \propto Id_{\C}$ (e.g. in one dimension), the terms for $k > 0$ vanish in the above equation and we get the simple well-known expression $X_\tau = (X_0 \delta_0 + I_\tau + dB_\tau )\ast h$.

% %\vspace{0.2cm}
\item  \textbf{Exponentionaly distributed delays:} Let us now treat the case $g: x \mapsto \beta e^{-\beta x} H(x)$. In this case, $\hat{g}(\xi) = \frac{\beta}{\beta + 2 i\pi \xi}$. Therefore, $\frac{2i\pi \xi}{\hat{g}(\xi)} = -2 \pi (\frac{2\pi}{\beta} \xi^2 -i \xi)$ which corresponds to the red curve in the left picture of figure \ref{fig: gaboules pictures}. Operators $A$ satisfying the spectral condition~\ref{eq: spectral condition} are the ones whose spectrum is contained in an open ball centered at $-l$ that does not intersect the red curve (blue disks of in figure \ref{fig: gaboules pictures}).

{When $A$ is negative, the operators $\U$ and $\V$ can be made completely explicit. Indeed, observing} that:\\
$\frac{\hat{g}(\xi)}{l \hat{g}(\xi) + 2i\pi \xi} = \frac{\beta}{(\frac{\beta}{2} + 2i\pi \xi)^2 - \frac{\beta^2\Delta^2}{4}} = \frac{\beta}{(\beta\frac{1 + \Delta}{2} + 2i\pi \xi) (\beta\frac{1- \Delta}{2} + 2i\pi \xi)}$
with $\Delta = \sqrt{1 - 4l/\beta}$, the operator $\V$ is the convolution operators generated by $\beta \big(h_- \ast h_+\big)$ with $h_\pm:t \mapsto e^{-{\beta\frac{1 \pm \Delta}{2}}t}H(t)$. Similarly $\U$ is generated by $\beta \big(h_- \ast h_+ + \frac{1}{\beta}h'_-\ast h_+\big)$. Even more explicitly, for $\beta > 4l$, $\V$ is generated by $t \mapsto \frac{2}{\Delta}e^{-\frac{\beta}{2}t}sh\big(\frac{\beta \Delta}{2}t \big)$ and $\U$ by $t \mapsto \frac{1}{\Delta}e^{-\frac{\beta}{2}t}\Big(\ sh\big(\frac{\beta \Delta}{2}t \big) + \Delta ch\big(\frac{\beta \Delta}{2}t \big) \Big)$. When $4l>\beta$, a similar result holds replacing the hyperbolic functions $ch$ and $sh$ by $\cos$ and $\sin$.

\item  \textbf{Single fixed delay:} For $g = \delta_0 + \alpha \delta_\theta$, we have $\hat{g}(\xi) = 1 + \alpha e^{2i\pi \theta \xi}$. The convergence domain of the expansion (condition~\ref{eq: spectral condition}) is shown in the middle and right pictures of figure \ref{fig: gaboules pictures}, for two different $\alpha \in \R_+$. The red curve seems to be living on the 2-dimensional projection of a simple 3-dimensional cone of revolution whose section is a circle. In that case, it appears quite difficult to express $\U$ and $\V$ in a simple form, though their Fourier transform is explicit.
\begin{figure}
\centering
\includegraphics[width=0.99\textwidth]{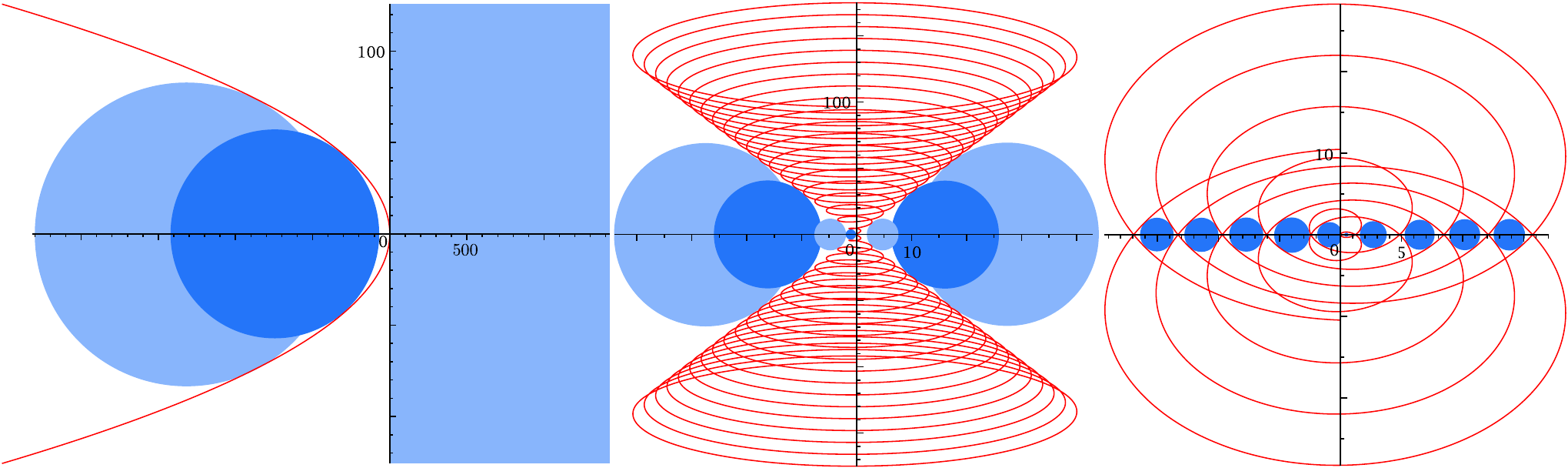}
% gaboules_pict.pdf: 587x175 pixel, 72dpi, 20.71x6.17 cm, bb=
\caption{ The three different pictures correspond to different time-convolution kernels $g$. The red curves are the parametric plots of $\xi \in \R \mapsto \frac{2i\pi \xi}{\hat{g}(\xi)} \in \mathbb{C}$ and the blue balls are examples of sets within which the eigenvalues of the space operator $A$ need to live for the expansion to be well-defined. The eigenvalues have to be contained in a single ball. The center of each ball is $-l$, for different $l \in \R$. To satisfy the spectral condition (\ref{eq: spectral condition}) the balls cannot intersect the red lines.
(left) Exponentialy distributed delays with $\beta = 2 \pi$ and $\xi \in [-20,20]$. (middle) Single delay with $\alpha = 2$ and $\theta = 1$ and $\xi \in [-20,20]$ (right) Single delay with $\alpha = 0.3$, $\theta = 1$ and $\xi \in [-5,5]$}
\label{fig: gaboules pictures}
\end{figure}
\end{itemize}
{{\bf Remark:} As illustrated in the previous example, a procedure to find the constant $l$ such that the expansion converges consists in plotting on the same figure the complex eigenvalues of the spatial operator and the red curve $\xi \in \R \mapsto \frac{2i\pi \xi}{\hat{g}(\xi)} \in \mathbb{C}$ related to the time operators. If there exists a ball centered on the real axis which contains all the eigenvalues and that does not intersect the red curve, then choosing $-l$ as the value of its center ensures that the expansion will converge.}
%\vspace{-0.4cm}
\subsection{{Stochastic heat equation}}
\label{sec:heat equation}
%\vspace{-0.4cm}
{Let us now deal with a classical the stochastic heat equation on $\mathbb{S}^1$ (described as the interval $[0,1]$ where $0$ and $1$ are identified) as a classical example of linear partial differential equations:
%\vspace{-0.15cm}
$$
\frac{\partial u}{\partial t}(x,t) = \Delta u(x,t) + v(x,t) + \sigma \eta(x,t)
$$}

%\vspace{-0.4cm}
\noindent {with periodic boundary conditions, where $\Delta$ is the Laplacian on $\mathbb{S}^1$, $v$ is an external forcing and $\eta$ is a multidimensional white noise. The input $v(x,t)$ is set to $\delta_{x=x_0}(x)$ and we take the initial condition $u(.,t=0)=0$.}

\begin{figure}
 \centering
 \includegraphics[width=0.7\textwidth]{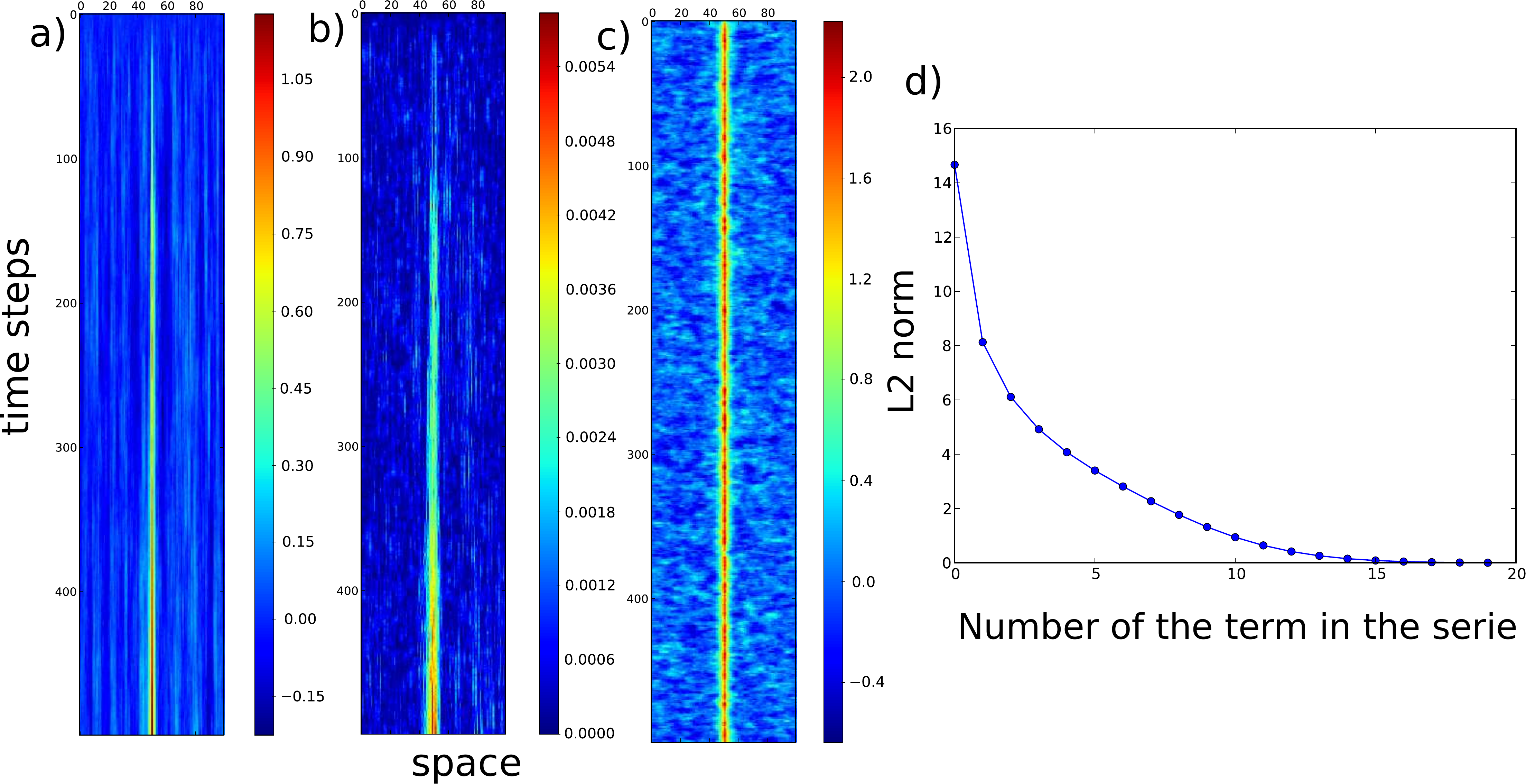}
 % heat_simus.pdf: 551x220 pixel, 72dpi, 19.44x7.76 cm, bb=0 0 551 220
 \caption{{Application of the expansion method to the stochastic heat equation on the circle with a Dirac source on the neuron in the middle. a) Space-time diagram of the solution given by the expansion method for $dt = 0.01$. b) Space-time diagram of the error between the solution in a) and the solution given by Euler's method. c) Space-time diagram of the solution given by the expansion method for $dt = 1$. d) $L_2$ norm of the terms in the expansion. The parameters for these simulations are $n = 100$, number of time steps $=500$, $\sigma = 0.1$ and $l=2$.}}
 \label{fig: heat simus}
\end{figure}

{The Laplacian operator has eigenvalues $-4\pi^2 k^2$ with $k\in \N$, corresponding to the eigenvectors $\cos(2k\pi x)$ and $\sin(2k\pi x)$. Since the eigenvalues are not bounded, it is not possible to find a suitable constant $l$ to define the solution of the heat equation in our framework. However, the semi-discretized in space equation overcomes this problem by preventing the existence of very fast oscillations (corresponding to large eigenvalues of the Laplacian). We choose to discretize the space with $N$ points regularly spaced, corresponding to a discretization step $dx=1/N$. The resulting equation corresponds to~\eqref{eq: differential SDE} in dimension $N$, with $g = \delta(t)$ and $A \in \R^{N \times N}$ such that $A_{ii} = -2/dx^2$, $A_{ij}=1/dx^2$ if $i=j \pm 1$, $A_{1n} = A_{n1} = 1/dx^2$ accounting for the periodicity of the medium, and $A _{ij} = 0$ otherwise . This matrix has eigenvalues in $[-4/dx^2,0]$. This suggests the choice $l=2/dx^2$ so that all the eigenvalues are in this ball a center $-l$ and radius $l$. This ball intersect the imaginary axis only in $0$ (corresponding to spatially constant functions), so convergence issues might arise if one of the terms $ W^k\cdot (\widetilde{I}_\tau + \Sigma\cdot dB_\tau)\cdot \U\cdot \V^k$ is spatially constant, which clearly never occurs in our case. Therefore, our expansion is well-posed and provides a numerical scheme to compute the solution as shown in figure \ref{fig: heat simus}.d. In that figure, we exhibit the fact that the solution is well retrieved by the expansion, and the error compared to Euler's scheme with a time step $dt=0.01$ (Fig. \ref{fig: heat simus}.b) is more than two order of magnitude smaller than the solution. An interesting point of this method is that it works for any time step interval $dt$ which is not the case for the Euler method which rapidly diverges as soon as the CFL condition is not satisfied for instance. Moreover, extending the approach to a delayed formalism $g \neq \delta$ is costless in our framework.}

%\vspace{-0.4cm}
\subsection{{Weakly connected networks}}
\label{sec:WCN}
%\vspace{-0.4cm}
{
Let us now go back to the problem of weakly connected neural networks introduced~\ref{sec: motivation}. The application of theorem \ref{thm: weakly expansion} to the Ornstein-Ulhenbeck process (see first example in section~\ref{sec:examples}), allows characterizing the equilibrium connectivity of system \eqref{eq: linear learning neural network} with a null initial condition as the solution of:
%\vspace{-0.1cm}
\[
W^* = \frac{1}{\kappa} \sum_{k,q=0}^{+\infty} \frac{{W^*}^k}{l^{k+1}}.I.\T_h^{k+1}.{\T_h^{q+1}}'.I'.\frac{{{W^*}'}^{q}}{l^{q+1}}
\]
%\vspace{-0.3cm}}

\noindent {In the weakly connected limit, i.e. $\kappa \to + \infty$, it can be shown that $\|W^*\| < l$ and we can find a good approximation of $W^*$ in the form
%\vspace{-0.3cm}
\[
W^* = \frac{I_m^2}{\kappa l^2} I.\T_h.\T_h'.I' + \mathcal{O}(\frac{I_m^2}{\kappa l^3})
\]}
%\vspace{-0.3cm}

\noindent{where $I_m = \sup_{t \in \R}  \|I(t)\|_2$. This expression makes explicit the fact that the equilibrium connectivity stores not only the spatial, but also the temporal correlation of the inputs. The method also makes it possible to go at arbitrary order and consider also correlated forcing noise together with delays in the communication term.
}

\medskip

\noindent{\bf Acknowledgements}
The authors would like to thank Francis Maisonneuve and Gilles Wainrib and also the anonymous reviewer who significantly hepled improving the paper.
MG and JT were partially funded by the ERC advanced grant NerVi \#
227747 and MG by the r\'egion PACA, France.
% Acknowledgements text here
%\vspace{-1.5cm}

% \begin{thebibliography}{00}
% % please try to use the bibitem system -
% % the references should be in alphabetical order of authors' names.
% % Articles with a single author first, author will 1 co-author next,
% % then author with several co-authors;
% 
% 
% % \bibitem{label}
% % Text of bibliographic item
% 
% \bibitem{label}
% 
% \end{thebibliography}

\end{document}